\documentclass[final,5p,times,twocolumn]{elsarticle}

\usepackage{amssymb}
\usepackage{amsmath}
\usepackage{float}

\usepackage{algorithm,algpseudocode}

\usepackage{xcolor}
\usepackage{listings}
\newcommand{\ii}{\mathrm{i}}
\newcommand{\dd}{\mathrm{d}}
\newcommand{\ham}{\mathcal{H}}
\newcommand{\e}[1]{\mathrm{e}^{#1}}

\newcommand{\ket}[1]{\ensuremath{\left|#1\right\rangle}}

\newcommand{\commie}[1]{\left[#1\right]_{-}}

\usepackage[hidelinks]{hyperref}

\newcounter{bla}

\makeatletter
\def\ps@pprintTitle{%
	\let\@oddhead\@empty
	\let\@evenhead\@empty
	\def\@oddfoot{}%
	\let\@evenfoot\@oddfoot}
\makeatother

\begin{document}
    
    \begin{frontmatter}

        \title{Parallel time integration using Batched BLAS (Basic Linear Algebra Subprograms) routines}

        \author[a]{Konstantin Herb\corref{author}}
        \author[a]{Pol Welter}
        \address[a]{Department of Physics, ETH Zürich, Otto-Stern-Weg 1, 8093 Zurich, Switzerland}
        
        \cortext[author] {Corresponding author.\\\textit{E-mail address:} science@rashbw.de}

        \begin{abstract}
            We present an approach for integrating the time evolution of quantum systems. We leverage the computation power of graphics processing units (GPUs) to perform the integration of all time steps in parallel. The performance boost is especially prominent for small to medium-sized quantum systems. The devised algorithm can largely be implemented using the recently-specified batched versions of the BLAS routines, and can therefore be easily ported to a variety of platforms. Our PARAllelized Matrix Exponentiation for Numerical Time evolution (PARAMENT) implementation runs on CUDA-enabled graphics processing units. 
            
        \end{abstract}
        
        
    \end{frontmatter}

    \section{Introduction}
    The last decade has seen the advent of quantum technologies. Significant advances have paved the way for promising applications in computing, sensing or communication. The time evolution of quantum systems is governed by the Schr\"odinger equation with a time-dependent Hamiltonian 
    \begin{equation}
        i \hbar \frac{d\psi}{dt} = H(t) \psi(t)
        \label{eq:schroedinger}.
    \end{equation}
    
    Efficiently solving eq.~\eqref{eq:schroedinger} is key to understanding the systems at hand, and to their successful technological application.
    Substantial effort has been put into simulating large quantum systems (thousands of degrees of freedom), by exploiting specific properties that allow a reduction of the Hilbert space. In addition, significant effort has been put on porting such simulations to GPUs. 
    
    Nevertheless, a lot of research focusses on small to medium-sized quantum systems. Although a simulation of such systems can be easily done on modern CPUs, computational expense scales with the number of time steps. As an example, consider an electron spin simulated in the laboratory (non-rotating) frame. Here, GHz control fields are applied for multiple microseconds leading to the need for integrating over tens of thousands of time steps. Often, this problem can be tackled by a suitable approximation (e.g. the rotating frame \cite{slichter1990principles}). However, such an approximation may not always be convenient or even possible: for instance, when studying non-secular effects like the Bloch-Siegert shift, the strong driving regime, complex modulated waveforms, or when performing a validation of the control software of an experiment.
    Accurate simulation rather than approximations (with sometimes elusive side effects) builds extra confidence in the correctness.
    Efficiency is key, too: when developing new pulse sequences for quantum control, it is desirable to run simulations at interactive speeds, such that the designer can quickly iterate different parameters. Furthermore, applications in the field of quantum optimal control rely on the fast evaluation of a time-dependent Hamiltonian.

    With the recent advent of artificial intelligence research, new powerful tools have emerged. For suitable tasks, a modern GPU offers the computational power of a supercomputer from the mid 2000s, off-the-shelf and at a fraction of the cost. Like supercomputers, GPUs focus on massive parallelization. The implementation of matrix-matrix multiplication in suitable hardware structures, combined with fast memory access, allows for a fast and parallelized tackling of a variety of computational tasks. 
    
    We present an approach for solving the time-dependent Schr\"odinger equation in a form that is frequently encountered in experimental realizations of a variety of physical quantum systems. Our approach leverages the parallelization of GPUs for integrating the steps of the time propagation in parallel and relies on a fast memory connection. We use the recently standardized Batched BLAS routines \cite{Dongarra2017} and a minimal set of custom functions. Therefore, our approach can be ported to a variety of platforms including, e.g., GPUs of other vendors or Field Programmable Gate Arrays (FPGAs). %
    
    This paper is structured as follows: First, we describe the concept of slice-wise propagation and the parallelization of the calculation. Then, we focus on the two underlying problems to be solved. This is followed by details about the implementation using the Batched BLAS functions. Furthermore, we improve the convergence order by extending the approach to a Magnus integrator. Finally, we showcase the runtime and the convergence using a suitable example of a driven two-level system.
    
    \section{Computational approach}
    For a system with a Hamiltonian $H = \hbar \mathcal{H}$ that is stationary in time, the exact solution of the Schr\"odinger or the von Neumann equation is given by 
    \begin{subequations}
        \begin{equation}
            \ket{\psi(t)} = \e{-\ii H t / \hbar} \ket{\psi(t=0)} = U(t) \ket{\psi(t=0)} \quad\text{and}
        \end{equation}
        \begin{equation}
            \rho(t) = U(t) \, \rho_0 \, U^\dagger(t) \text{,}
        \end{equation}
    \end{subequations}
    respectively. Here, $\ket{\psi}$ denotes the system state of a pure quantum system and $\rho$ the density matrix of a mixed state.  One approach to treat arbitrary time-dependent Hamiltonians is to`slice' the Hamiltonian, known also as Euler's method. During the finite duration $\Delta t$ of each slice, the Hamiltonian is assumed to be stationary. For sufficiently small equidistant time steps, the sequence
    \begin{equation}
        U(t) = \lim\limits_{\Delta t \rightarrow 0} U^{(n)} \cdot U^{(n-1)} \cdot ... \cdot U^{(0)} = U(t_n)\mathrm{.}
        \label{eq:slices}
    \end{equation}
    with
    \begin{equation}
        U^{(n)} = \exp(-i\mathcal{H}^{(n)}\Delta t)
    \end{equation}
    converges to the true solution of the original problem. If the dimension of the Hamiltonian is large, Krylov methods can be employed to simplify the exponentiation. Here, we instead focus on small, dense matrices.
    
    The algorithm thus involves two steps: (i) calculating the matrix exponential efficiently for many matrices and (ii) executing the matrix multiplication of the individual slice-propagators. In the following subsections, we will discuss the chosen approaches for each of the two aspects.
    
    We assume that the time-dependent Hamiltonian takes the form
    \begin{equation}
        \mathcal{H}(t) = \mathcal{H}_0 + \sum_{i=1}^{N} c_i(t)\,\mathcal{H}_i
        \label{eq:format_control}
    \end{equation}
    where $ \mathcal{H}_0$ is typically named drift Hamiltonian and $\mathcal{H}_i$ and $c_i(t)$ denote control Hamiltonians and control coefficient array. The control terms represent the coupling of the quantum system, e.g., time-varying magnetic or electric fields (control fields). The decomposition in eq.~\ref{eq:format_control} does not trade any generality: $\mathcal{H}_i$ are simply a basis of the $\mathcal{H}(t)-\mathcal{H}_0$ space. Numerous experimental situations however conform very well to the form of eq.~\eqref{eq:format_control}, with a number $N$ of control fields that is much smaller than the total dimension of the Hamiltonian space.
    
    \begin{figure}[h]
        \includegraphics[width=\linewidth]{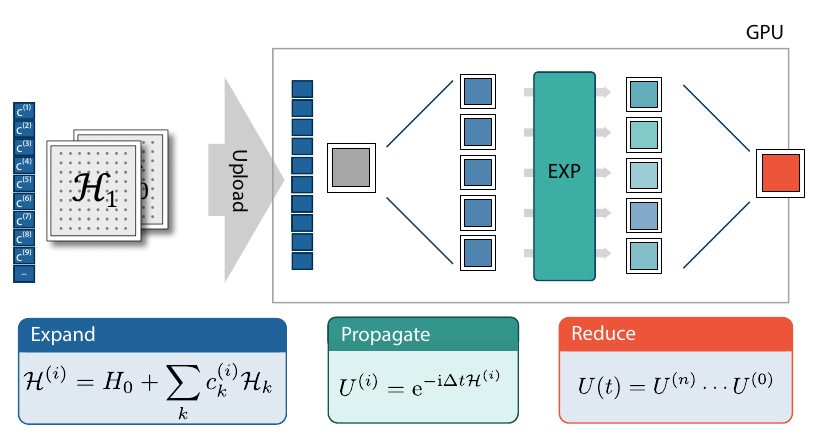}
        \caption{Parallel time integration of the PARAMENT integrator. We upload a minimal amount of data to the GPU and propagate all time steps synchronously. The final propagator is obtained by iterative pair-wise multiplication of the slice propagators. In the case of the Magnus implementation, the expansion step also calculates the coefficients for the commutator terms.}
        \label{fig:scheme}
    \end{figure}
    The overall structure of the approach is shown in Figure \ref{fig:scheme}: We transfer a minimal amount of data to the GPU, the control arrays and the Hamiltonians. Then we expand the Hamiltonians for each time step in the GPU memory. Subsequently we exponentiate all slice Hamiltonians in parallel. For the final propagator,  we reduce the slice propagators by repeated pair-wise matrix multiplications.  
    
    \subsection{Matrix exponential}
    The numerical computation of the matrix exponential has been treated extensively by Moler and Van Loan in their famous `19 dubious ways' paper \cite{Moler2003}. Moler and Van Loan describe six main classes of algorithms: (1) Series-based methods relying on, e.g., Taylor or Pad\'e approximations, (2) Methods relying on ODE solvers, (3) polynomial methods which are typically not very attractive due to their high computational cost, (4) matrix decomposition methods, (5) splitting methods like the powerful scaling-and-squaring technique and (6) Krylov methods, which are interesting for large matrices. On CPUs, matrix decomposition methods are particularly powerful due to well-established implementations like the Schur decomposition in LAPACK using the ZGEES or CGEES function \cite{laug}. This approach is especially competitive for Hermitian matrices like Hamiltonians, where the decomposed matrix is always diagonal. It is used e.g. in the MATLAB software package. In general, for non-Hermitian matrices, most modern implementations (e.g. in Python and MATLAB) combine the scaling and squaring technique with a series method, either Taylor or Pad\'e. Typically the matrix is first scaled down until its norm is sufficiently small so that its exponential can be well approximated by a (reasonably) truncated Taylor or Pad\'e approximation. The exponential of the original matrix can be recovered by raising the small-norm exponential to the correct power.
    
    For our implementation of the highly-parallelized calculation of the matrix exponential we choose a different series approach based on the Chebyshev polynomials. The expansion of the matrix exponential in a Chebyshev series in the context of time propagation of quantum systems goes back to Tal-Ezer and Kosloff \cite{TalEzer1984}. The expansion of the complex exponential in a Chebyshev series can be obtained by starting with the expansion of $\e{\ii \omega x}$ on the unit interval $x \in [-1,1]$ 
    \begin{equation}
        \e{\ii \omega x} = a_0 + 2 \sum_{k=1}^{\infty} a_k T_k(x) \qquad \text{with} \qquad a_k = \ii^k J_k(\omega) 
        \label{eq:chebyshev_scalar}
    \end{equation} where the coefficients are readily obtained from Abramovic Stegun (9.1.21) and $T_k(x)$ and $J_k(\omega)$ denote the Chebyshev polynomials and the Bessel functions of the First Kind respectively. An approximation of $\e{-\ii \zeta}$ for the interval $\zeta \in [\alpha,\beta]$ can be obtained by a straightforward affine-linear transformation of $x$. The scalar expansion \eqref{eq:chebyshev_scalar} can be extended to matrix arguments and reads in our case for the slice propagators \cite{Lubich2008}
    \begin{equation}
        \e{-\ii G} \approx
        \e{-\ii(\alpha+\beta)/2} \left[a_0 \mathbb{I} + 2 \sum_{k=1}^{m_\mathrm{max}} a_k T_k\left(\frac{2}{\beta-\alpha}\left(G-\frac{\alpha+\beta}{2}\mathbb{I}\right)\right)\right] \text{.}
        \label{eq:chebychev_series}
    \end{equation}
    Here, $G=\ham \Delta t$ is Hermitian and we suppose that the spectrum of this matrix lies between $\alpha = \lambda_\mathrm{min}(\ham\Delta t)$ and $\beta = \lambda_\mathrm{max}(\ham\Delta t)$. For the transformed expansion, the coefficients $a_k$ read
    \begin{equation}
        a_k = (-\ii)^k J_k\left(\frac{\beta-\alpha}{2}\right) \text{.}	
    \end{equation} 
    For the numerical evaluation of the matrix exponential, we will truncate the series in \eqref{eq:chebychev_series} at $m_\mathrm{max}$. Lubich \cite{Lubich2008} estimated the error for approximating the matrix exponential by using a Chebyshev series truncated after the $m$'th term to
    \begin{equation}
        \epsilon = 4 \left(\e{1-\left(\frac{\beta-\alpha}{4m+4}\right)^2}\frac{\beta-\alpha}{4m+4}\right)^{m+1} \text{.}
        \label{eq:error_chebychev}
    \end{equation}
    We use \eqref{eq:error_chebychev} to determine $m_\mathrm{max}$ under the requirement that $\epsilon$ is smaller than the machine precision of the respective datatype. This gives us constrains for the maximum norm that the argument can take, c.f. Table \ref{tab:norm_boundries}.
    
    \begin{table*}
        \resizebox{\textwidth}{!}{
            \begin{tabular}{l|llllllllllll}
                $m_\mathrm{max}$ & 3 & 5 & 7 & 9 & 11 & 13 & 15 & 17 & 19 & 21 & 23 & 25 \\
                \hline
                $||G||_\mathrm{max}$ FP32 & $0.033$ & $0.219$ & $0.620$ & $1.218$ & $1.980$ & $2.873$ & $3.873$ & $4.959$ & $6.118$ & $7.336$ & $8.606$ & $9.919$	 \\
                $||G||_\mathrm{max}$ FP64 & $2\times 10^{-4}$ & $0.008$ & $0.050$ & $0.163$ & $0.368$ & $0.677$ & $1.088$ & $1.596$ & $2.194$ & $2.874$ & $3.627$ & $4.447$ \\			
            \end{tabular}
        }
        \caption{Maximum norms of the exponent of a matrix exponential that guarantees machine precision when calculating the matrix exponential by using the Chebyshev expansion \eqref{eq:chebychev_series} truncated at $m_\mathrm{max}$. Here, we approximated the spectral range of $G$ with $\alpha = - ||G||$ and $\beta = ||G||$ and solved \eqref{eq:error_chebychev} for $||G||$. For calculations in single-precision floating-point format (FP32) we required  $\epsilon < 2^{-24}$, for the double-precision floating-point format (FP64) $\epsilon < 2^{-53}$.}
        \label{tab:norm_boundries}
    \end{table*}
    \begin{algorithm}[H] 
        \caption{Matrix exponential}
        \label{alg:loop}
        \begin{algorithmic}[1]
            \Require{Hamiltonian $\ham$, time step $\Delta t$, $\qquad\qquad\qquad\quad$ spectral boundaries $\alpha,\beta$}
            \Ensure{Propagator $U=\e{-i\Delta t \ham}$}
            \Statex
            \State{$X = \frac{2}{\beta-\alpha}$ $\left(\Delta t \ham - \frac{\alpha+\beta}{2}\mathbb{I}\right)$}
            \State{$D_{m_{max}+2} = 0$}
            \State{$D_{m_{max}+1} = 0$}
            \For{$k = m_{max}$ downto $0$}                    
            \State {$a_k = (-\ii)^k J_k\left(\frac{\beta-\alpha}{2}\right)$}
            \State {$D_k = a_k \mathbb{I} + 2 X D_{k+1}-D_{k+2}$}
            \EndFor
            \State \Return {$\e{-\ii (\alpha+\beta)/2} (D_0 - D_2)$}
        \end{algorithmic}
    \end{algorithm}
    The series in \eqref{eq:chebychev_series} can be efficiently computed by the Clenshaw algorithm \cite{Clenshaw1955}.
    
    The number of necessary steps (and thus the number of matrix multiplications) grows at least linearly with the matrix norm. For our application, this is not a major problem as the matrix norm per time slice $||\mathcal{H}\Delta t||$ is expected to be small if the time step is sufficiently short. In a general context, Chebyshev-series-based matrix exponentiation has been successfully combined with the powerful scaling-and-squaring approach \cite{Auckenthaler2010}. This approach brings down the number of necessary matrix multiplications to a growth that is only logarithmic in the matrix norm. It is worth noting that the Chebyshev approach is also very appealing for problems with sparse Hamiltonians, as it relies only on matrix products of the form $\mathrm{sparse} \times \mathrm{dense}$. This circumstance could be further leveraged for bigger Hamiltonians. Nonetheless, our approach requires \textit{a priori} knowledge about the expected spectrum of the Hamiltonian.
    As the boundaries $\alpha$ and $\beta$ of the Eigenvalue range are often not available, we estimate them by using a suitable matrix norm (Gershgorin's theorem).

    \subsection{Reduction by matrix multiplication}
    After obtaining all slice propagators, we multiply them together to obtain the full evolution, c.f.~eq.~\eqref{eq:slices}. We can again parallelize this step by leveraging the associativity of the matrix multiplication. Our slice propagators $U^{(i)}$ are all quadratic and have the same dimensions. Therefore, we simply compute pair-wise products of consecutive slice propagators: $U^{(0)}\cdot U^{(1)}, U^{(2)}\cdot U^{(3)}, ..., U^{(n-1)}\cdot U^{(n)}$. By repeating this scheme $\log_2(n)$ times, where we halve the number of remaining time slices with every round, we compute the overall product. As we will see in the next section, the pair-wise matrix multiplications can be easily parallelized using the Batched BLAS functions. Although there exist more sophisticated approaches to this problem in the literature \cite{Ladner1980,Irony2004}, during our numerical tests it turned out that this simple approach is sufficient as the major computational cost is given by the exponentiation step.
    
    \section{Implementation}
    The main idea of our implementation is that the computationally expensive exponentiation performed by algorithm \ref{alg:loop} is run in parallel for all time steps. As all time steps have approximately a similar spectrum we can choose global $\alpha$, $\beta$ and $m_{max}$ that are suitable across all time steps. Together with the fact that all matrices have the same shape, this makes the algorithm ideally suitable for a Single Instruction Multiple Data (SIMD) platform, such as a GPU. Furthermore, algorithm \ref{alg:loop} has been designed in a way that it can be implemented by only a small subset of the BLAS functions (batched and non-batched) which will bear the majority of the computational workload.
    
    \subsection{Batched BLAS routines}
    The Batched BLAS routines have been standardized in 2017 \cite{Dongarra2017} and were created with the intend to maintain high compute resource utilization when dealing with small-to-medium-sized matrices. Looking at the ubiquitous General Matrix Multiply (GEMM) routine,  the batched version performs the operation
    \begin{equation}
        C[k] = \alpha\,\, A[k]\cdot B[k] + \beta \, C[k] \qquad \forall k
    \end{equation}
    for a batch of length $k$ in parallel. 
    
    It is obvious that the Batched BLAS functions can be used to parallelize the matrix exponentiation step. However, the flexibility of the GEMM routine allows us to reuse the GEMM function for nearly all steps in Figure \ref{fig:scheme}. As $A$ and $B$ are pointers, the strided batched GEMM routine can also be reused for the pair-wise reduction operation, where we double the memory stride to twice the propagator size. $A$ and $B$ can point to the same region in memory with only an offset of one matrix between them. Depending on the implementation, special attention may have to be brought to possible bottlenecks induced by memory miss-alignments. However, our tests with the NVIDIA cuBLAS implementation (see below) did not reveal major negative effects.

    Similarly, the addition of the Bessel coefficients to the diagonal elements of the propagators during the exponentiation step (lines 5 and 6 in algorithm \ref{alg:loop}) can be implemented using the batched version of AXPY (scalar $\times$ vector multiplication, "$a\cdot x$~plus~$y$"). By this, the full integration can be done solely with BLAS routines, without any custom compute kernels required. This makes the approach easily portable to a variety of systems.

    \subsection{GPU implementation using NVIDIA cuBLAS}
    We implemented the proposed integration scheme using the CUDA platform and by using the cuBLAS library. The batched version of the AXPY routine had to be implemented by a custom CUDA kernel as it is not yet part of the cuBLAS library (as opposed to other GPU frameworks like AMD's rocm). To reduce the memory requirements for working arrays, we implemented two Chebyshev iterations per loop iteration in algorithm~\ref{alg:loop}, c.f. Appendix.
    
    The performance of the integrator benefits from the fast memory connection on the GPU and the sheer number of compute processors available.
    On a Quadro P2000 GPU, for 80'000 time steps and a $12\times 12$ system, the exponentiation step is by far the most time consuming ($\sim 80\%$ of the total run time). The NVIDIA visual profiler (NVVP) reports that we make good use of the available resources (a memory bus utilization of $55\%$, a compute resource utilization of $85\%$, and an occupancy of $90\%$). During the subsequent reduction step, these numbers are slightly reduced  ($55\%$, $75\%$, and $35\%$). This may indicate room for optimization, but any improvements here are unlikely to improve total runtime significantly, since time expended for the reduction step is short already.
    It is worth noting that memory misalignment problems can decrease the efficiency of the Batched BLAS functions when e.g. setting the memory stride to unusual values. However, for the cuBLAS library we do not observe this to be a major problem.
    
    The user provides the control amplitudes as an array sampled at equidistant time points. From the control amplitude array we compute the actual exponents $-\ii G$, optionally by averaging the control vector over three points according to Simpson's quadrature rule. This makes sense in the context of a higher-order Magnus integrator (see section~\ref{sec:magnus} below). Otherwise, we simply set $G= \mathcal{H}\Delta t$ with $\mathcal{H}$ as per equation (\ref{eq:format_control}).
    
    The exact computation of the spectral boundaries $\alpha$ and $\beta$ can be expensive. We instead use the operator norm $||G||_1$ as an inexpensive upper bound. We use a single value for $\alpha$ and $\beta$ for all time steps, so we further bound the norm by the triangle inequality: $||G||_1/\Delta t < ||\mathcal{H}_0||_1 + \sum_{i=1}^{N} |c_i(t)|\,||\mathcal{H}_i||_1$. By requiring the user to keep the control amplitudes $c_i(t) \in [-1,1]$, rescaling the control Hamiltonians if necessary, we hence define
    \begin{equation}
        - \alpha = \beta =  \Delta t  \sum_{i=0}^{N}||\ham_i||_1 \text{.}
    \end{equation}

    \section{Magnus integrators}
    \label{sec:magnus}
    It can be shown that the error of the above integrator in a single time step scales as $\mathcal{O}(\Delta t^3)$, making it a second order integrator over the full integration time $T=N\Delta t$.
    The order of the convergence can be improved significantly by Magnus expansion \cite{Magnus1954, Blanes2009}. Magnus integrators have been successfully used for solving the Schr\"odinger equation in various implementations, including on GPUs \cite{Auer2018}. However, these typically again focus on large Hamiltonians, rather than the small to medium-sized problems with long time horizon that we are considering.
    
    Magnus proposed that the solution $U(t)$ to any time-dependent ODE (like the Schr\"oedinger equation) can indeed be written as the exponential of some matrix $\Omega(t)$. He provided a series expansion of that exponent, i.e. 
    \begin{subequations}
        \label{eq:magnus}
        \begin{align}
            U(t) &= \e{\Omega(t)} \quad \text{with} \quad \Omega(t) = \sum_i\Omega_i(t) \tag{\ref{eq:magnus}} \\
            \Omega_1(t) &= -\ii \int_0^t \dd t_1 \, \ham(t_1)  \text{,}	 \\
            \Omega_2(t) &= -\frac{1}{2}\int_0^t \dd t_1 \int_0^{t_1} \dd t_2 \, [\ham(t_1),\ham(t_2)]_{-} \text{,}\\ 
            \dots & \tag*{}
        \end{align}
    \end{subequations}
    where $[A,B]_{-} :=AB-BA$ is the commutator of the matrices $A$ and $B$. The interpretation is obvious: The term $\Omega_1(t)$ averages the Hamiltonian over the time slice. By replacing the integral with the mid-point rule, we recover the naive approach where we assume the Hamiltonian to be constant over the whole time slice.
    The higher-order terms $\Omega_2(t)$, $\Omega_3(t)$, ...~describe the effects when the Hamiltonian at time $t_1$ does not commute with itself at a later time $t_2$.
    
    For constructing a fourth-order integrator, it is sufficient to only include terms up to $\Omega_2$, and then to approximate the integrals with quadrature rules of sufficiently high order \cite{Blanes2000434}. The inclusion of this extra term neatly integrates with the method described above. As we will show next, including $\Omega_2$ is equivalent to simply adding extra control terms in equation (\ref{eq:format_control}). The extra terms correspond exactly to the commutators of the existing control Hamiltonians.
    
    Assume we have sampled the Hamiltonian at 3 equidistant time points, $\ham^{(1)} = \ham(0)$,  $\ham^{(2)} = \ham(\Delta t)$ and  $\ham^{(3)} = \ham(2\Delta t)$. Following Blanes et al.\footnote{c.f. eq. (256) of \cite{Blanes2009}}, we approximate the integrals in equations (\ref{eq:magnus}a) and (\ref{eq:magnus}b) by Simpson's rule and the trapezoidal rule respectively. The new exponent then reads
    \begin{equation}
        \begin{split}
            \Omega &= -\ii G \\ &= -\ii \left(\frac{\Delta t}{3}\left(\ham^{(1)} + 4 \ham^{(2)} +\ham^{(3)}\right) - \ii \frac{\Delta t^2}{3} \commie{\ham^{(1)},\ham^{(3)}} \right).
        \end{split}
    \end{equation}
    Inserting equation \ref{eq:format_control}, we find
    \begin{equation}
        \begin{split}
            G = 2 \Delta t\,\ham_0 &+ \Delta t\sum_{k=1}^N \left(\frac{1}{3}c_k^{(1)}+\frac{4}{3}c_k^{(2)}+\frac{1}{3}c_k^{(3)}\right) \underbrace{\ham_k}_*  \\
            &-\ii \frac{\Delta t^2}{3} \sum_{k=1}^N \left(c_k^{(3)}-c_k^{(1)}\right) \underbrace{\commie{\ham_0,\ham_k}}_* \\
            &-\ii \frac{\Delta t^2}{3} \sum_{k<k'}\left(c_k^{(1)}c_{k'}^{(3)}-c_k^{(3)}c_{k'}^{(1)}\right) \underbrace{\commie{\ham_k,\ham_{k'}}}_*.
        \end{split}
        \label{eq:magnus_coeffs}
    \end{equation}
    The matrix $G$ is still of the same form as $\ham$ in eq.~(\ref{eq:format_control}). The matrices marked with an asterisk now correspond to the new effective control Hamiltonians. Also note that compared to eq.~(\ref{eq:format_control}), we have now doubled the time step to $2\Delta t$. This is because Simpson's rule requires an extra sample in the middle of the interval.
    
    By exponentiating $-\ii G$ instead of $-\ii \ham\Delta t$ we thus double the time step, yet also increase the order of convergence of the integrator from 2 to 4. The price to pay is that the transform introduces additional control Hamiltonians. It maps $N$ original control Hamiltonians (and their corresponding control amplitude array) to $3/2 N  + N^2/2$ effective control Hamiltonians. This can be costly if the initial number of control Hamiltonians is already large. 
    
    Note that the transform $\ham \rightarrow G$ is entirely contained in the `expand' step in figure \ref{fig:scheme}. The time expended here is usually negligible though, as the computational bottleneck is the subsequent exponentiation. In practice, we find that the improved convergence allows for much fewer timesteps (for a given target accuracy), such that \emph{both} the expansion and propagation steps are proportionally faster to compute.\footnote{In reality, doubling the timestep improves the run-time by slightly less than a factor 2. This is because the exponent norm scales with the time step, so more iterations in the Chebychev expansion are required (see Table \ref{tab:norm_boundries}).}
    
    Finally, an extension to higher-order Magnus terms is possible and integrators of 6\textsuperscript{th} or 8\textsuperscript{th} order can be obtained, although the number of necessary commutator terms grows very fast. Here, the recently developed commutator-free Magnus integrators \cite{Alvermann2011} may provide an alternative.
    
    \section{Numerical experiments}
\begin{figure}[h]
        \centering
        \includegraphics[width=1.05\linewidth]{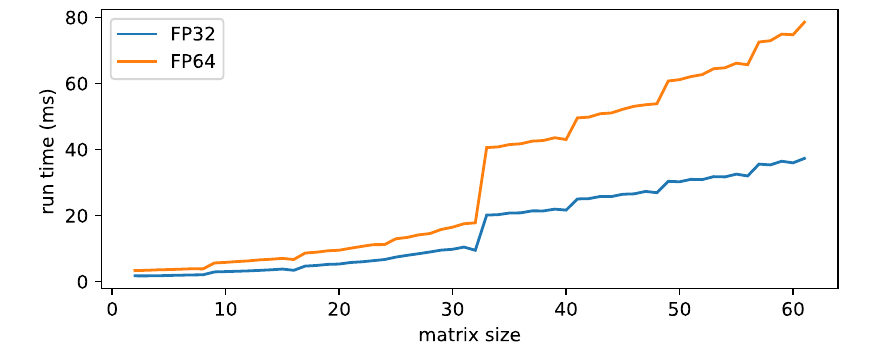}
        \caption{Runtime vs. matrix size for 80'000 time steps on the GPU. For small matrices ($d<32$), we see a very gentle increase in runtime when increasing the matrix size. At $d=32$ the cuBLAS library switches to a different compute kernel.}
        \label{fig:performance_size}
    \end{figure}

    We first test the performance by propagating several dense Hamiltonians as a function of matrix size. We used an NVIDIA V100 GPU. The results are shown in Figure \ref{fig:performance_size}.
    Next, we compare it against a fully parallelized CPU implementation, running on an Intel i9-9900 and linked against Intel MKL. The CPU implementation combines the expansion and propagation step in Fig. \ref{fig:scheme} into a single parallelized  function. This improved the performance due to more efficient caching as a direct port of the GPU code would be an unfair comparison. For a $12 \times 12$ Hamiltonian system (encountered often in nitrogen vacancy research) and 80'000 time steps, we find that the GPU runs 50 times faster for single precision (FP32) and 25 times faster for double precision (FP64).
    
\begin{figure}[h]
        \centering
        \includegraphics[width=1.07\linewidth]{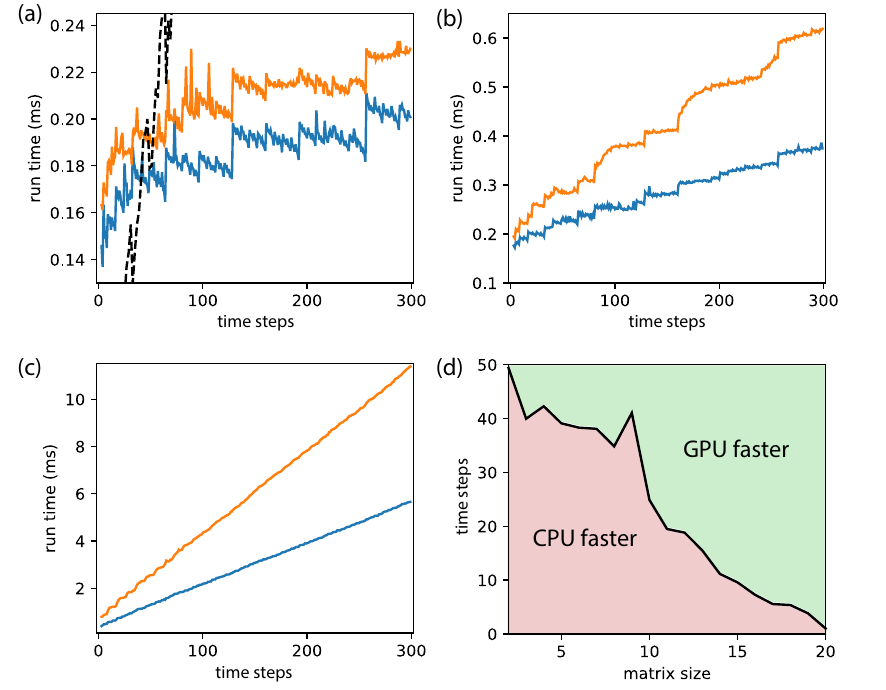}
        \caption{Perfomance as a function of time steps for various regimes of matrix dimensions. Blue traces are single precision (FP32) and orange traces are double precision (FP64). (a) shows the run-time for a $2 \times 2$ matrix. The black trace shows the measured run-time for a propagation implemented on a CPU (FP64). (b) shows the measured run-times for $40 \times 40$ matrices and (c) for $192 \times 192$. (d) shows the intersection of the run-time curves of the CPU and the GPU in the FP64 case.}
        \label{fig:performance_points}
    \end{figure}

    In the regime of a small number of time steps and a small matrix dimensionality, the CPU is faster as the "outsourcing" of the computation to the GPU comes with several overheads. This is highlighted in Figure \ref{fig:performance_points}. For small matrices, we see an approximately logarithmic increase in the runtime, due to the extra kernel launches required in the "reduce" step.
    
\begin{figure}[h]
        \includegraphics[width=\linewidth]{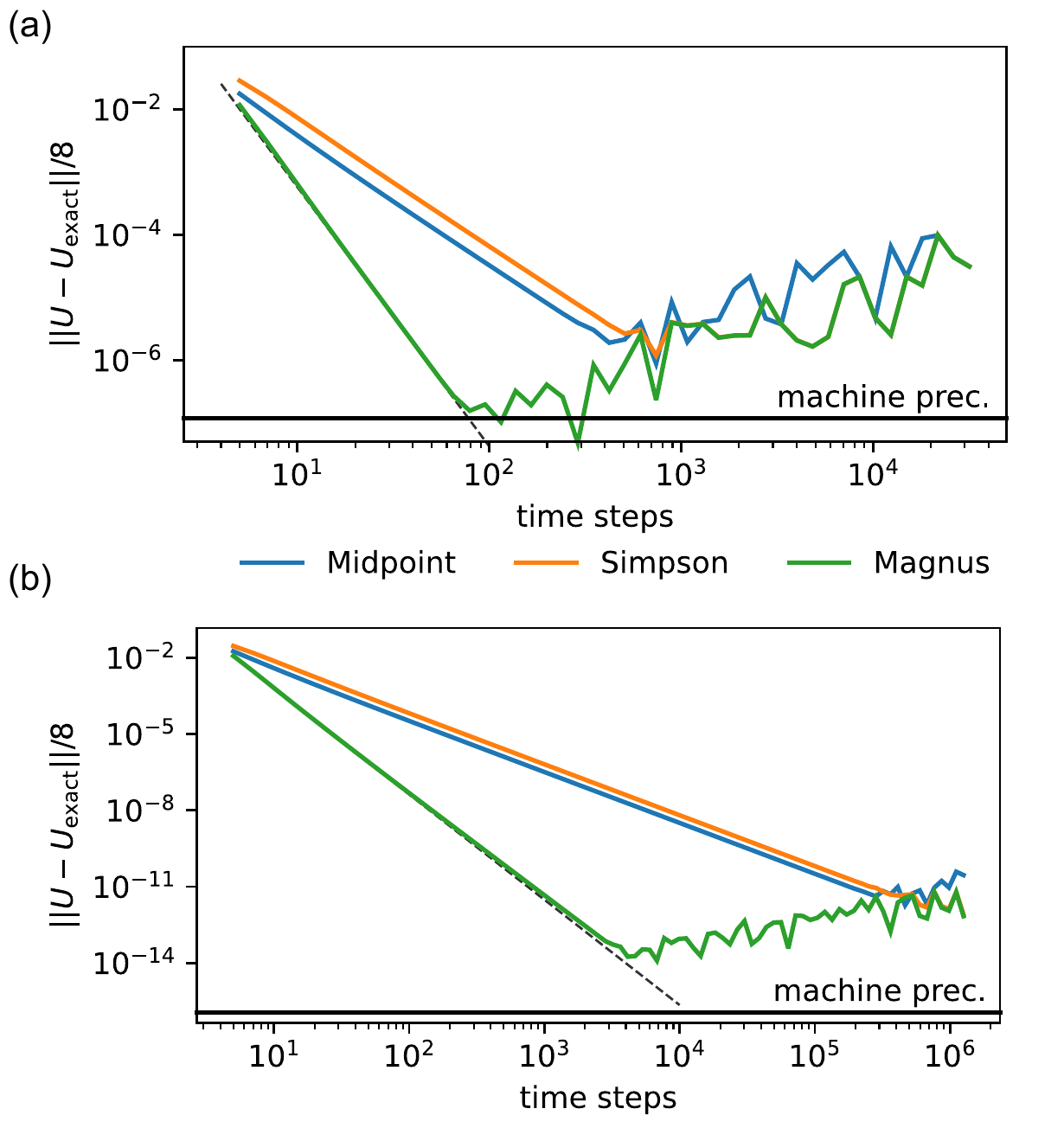}
        \caption{Convergence of the propagator for a qubit driven with a circularly polarized excitation field for various numbers of time steps (a) in case of single precision (FP32) accumulation and (b) double precision (FP64) accumulation. The parameters in \eqref{eq:driven_hamiltonian} are $\omega_0 = 1.0, \omega_1 = 0.1, \omega_\mathrm{rf} = 1.0$ and we evaluated the propagator at final time $t = 6.0$. The manually added black dashed-line indicates the same convergence order for the FP32 and the FP64 case}
        \label{fig:convergence}
    \end{figure}

    Next, we test the rate of convergence. We used the model of a qubit driven with a circularly polarized excitation field, a problem for which an exact analytical solution exists. The studied system Hamiltonian is 
    \begin{equation}
        \ham (t) = \frac{\omega_0}{2} \sigma_z + \cos(\omega_\mathrm{rf} t) \frac{\omega_1}{2}\sigma_x  + \sin(\omega_\mathrm{rf} t) \frac{\omega_1}{2}\sigma_y \text{,}
        \label{eq:driven_hamiltonian}
    \end{equation}
    where the $\sigma_i$ denote the Pauli matrices. Figure \ref{fig:convergence} shows the resulting accuracy when increasing the number of time steps while keeping the total evolution time fixed. We clearly observe that for both, FP32 and FP64, we achieve a better convergence when implementing a Magnus integration scheme while the computational cost per time step only increases marginally. For the sinusoidal drive, our implementation reaches machine precision for $\approx10^2$ time steps per oscillation cycle in the case of FP32, and for $\approx10^4$ time steps per cycle in the FP64 case. Furthermore, we see that when truncating the Magnus expansion at $\Omega_1$, increasing the quadrature order alone does not improve the convergence. 
    
    We observe that the error reaches a minimum after a certain number of timesteps. Finer timesteps do not improve the accuracy. Past this point, the error per time slice is limited by the machine precision, which accumulates over an increasing number of slices. With single-precision arithmetic, the achievable error is on the order of $10^{-5}$. This result might seem unsatisfying at first. However, this precision still by far exceeds the accuracy achieved in many experimental realizations of the quantum systems that the algorithm is designed to simulate. For instance, in the field of Nitrogen Vacancy (NV) center research, available Signal-to-Noise Ratio (SNR) frequently limits the practical \emph{experimental} accuracy to $\sim 1 \%$. In turn, this means that even consumer GPUs that throttle FP64 performance can provide acceptable performance.

    \section{The PARAMENT library}
    We provide the full-GPU integrator as a C library. We name the library PARAMENT (\textbf{PAR}allelized \textbf{M}atrix \textbf{E}xponentiation for \textbf{N}umerical \textbf{T}ime-evolution) and available for download\footnote{\url{https://github.com/parament-integrator/parament}}. The compiled DLL for the Windows platform and the UNIX version of the shared library can be included into a large range of applications including  Python, Julia, Matlab or LabVIEW.
    
    \begin{figure}[h]
        \includegraphics{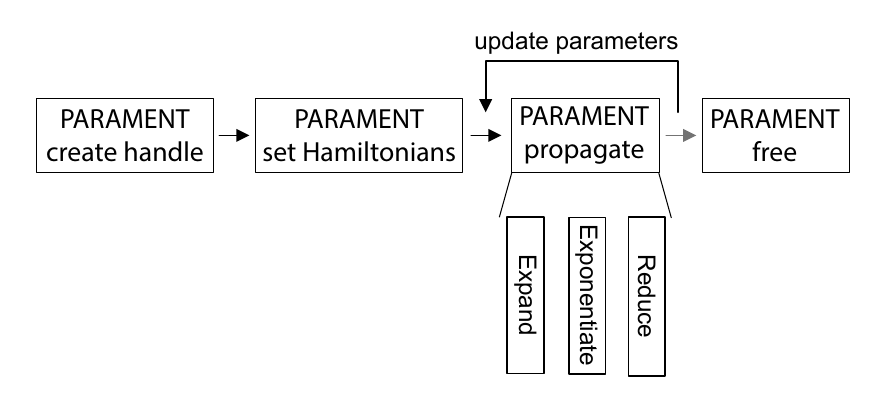}
        \caption{State-machine of the PARAMENT integrator.}
        \label{fig:state_machine}
    \end{figure}
    
    The usage model follows the steps of Figure~\ref{fig:state_machine}. First, the user initializes the integrator by calling the \texttt{Parament\_create()} function. It returns a handle to a newly created context. The context stores the state of the integrator. Multiple contexts may exist at any time, and used independently; however they are not thread-safe nor reentrant. The context must eventually be released by calling \texttt{Parament\_free()}.
    
    Next, the Hamiltonians are uploaded to the GPU using the \texttt{Parament\_setHamiltonian()} function. Here, the user also decides on the use of the Magnus expansion. When enabled, PARAMENT will then calculate the necessary commutators (effective control Hamiltonians) and upload them to the GPU as well.
    
    To obtain a propagator, the user calls the \texttt{Parament\_equiprop()} routine with the coefficient arrays. The Hamiltonians persist between propagations, so that repeated runs (e.g. with different control fields) are possible. Lastly, the implementation exposes several helper functions which allow tweaking the underlying numerics, e.g., the selection of a different $m_\mathrm{max}$. For a full description of available functions, the user is referred to the documentation of PARAMENT bundled together with the source code, or available on the project website.
    
    \subsection{Python}
    While written in C++, the PARAMENT library can easily be used together with other programming languages, including Matlab or LabView. For Python, we provide a reference binding, called \texttt{pyparament}. 
    This use-case was the initial motivation for the development of PARAMENT: A fast lab-frame propagation in interactive compute sessions, e.g., during the development of new microwave control schemes, possibly in Jupyter notebooks. Applications include, for instance, testing of advanced control schemes e.g. by frequency-modulating the microwave pulses \cite{Silveri2017}. It truly embraces the mindset of `interactive super-computing' \cite{Reuther2018}.
    
    Finally, \texttt{pyparament} is also compatible with the QuTIP framework\cite{Johansson2013}.
    
    The source code is available on the PARAMENT Github, or on the PyPI package repository.

    \section{Applications and outlook}
    The presented speed-up of lab-frame simulations will facilitate testing new control schemes during the sequence design of quantum control experiments. Due to implementation with only a few BLAS functions, our scheme can be ported easily to a variety of platforms, including AMD GPUs and FPGA devices.
    
    The applications of our integration scheme, however, go beyond what it was initially designed for. It is generally suitable for studying the evolution of small-to-medium-sized quantum system under complex-modulated control fields. It can be used to quickly determine the Floquet states and quasi energies of strongly periodically driven system, by diagonalizing the propagator \cite{Creffield2003}.
    
    As the size (degrees of freedom) of the simulated system increases, the presented approach requires vastly more memory. Fortunately, modern GPUs provide ample amounts thereof. If that is insufficient, PARAMENT could easily be adapted to scale across multiple GPUs. We hope to encourage the implementation of Batched BLAS routines that natively support multi-GPU calculations via fast GPU-to-GPU buses like the NVIDIA NVLink.
    
    A more visionary application could be the verification of an experimental control software and hardware stack. The quantum system could be replaced by PARAMENT, a fast digitizer, and a fast arbitrary waveform generator (AWG). If the latter two have direct access to GPU memory, performance may be sufficient to fully emulate the quantum system in near-real time.
    
    Our approach is not limited to quantum systems. It can be applied to any differential equation that can be cast into the form of equation \eqref{eq:schroedinger}. Of course, computational efficiency is best if the problem can be formulated with only a few `control Hamiltonians', such that limited upload rate from host computer to GPU does not affect the overall computation time.
    
    Finally, the fully GPU-based integration algorithm can be used as an essential building block for fully GPU-based optimal control optimization routines like the GRAPE algorithm \cite{Khaneja2005}. The optimizer must evaluate the time-evolution operator in every step of the optimization routine. We are confident that GRAPE can be greatly accelerated if built around PARAMENT. Depending on how large the norms of the slice Hamiltonians are, here the introduction of a scaling and squaring approach might be needed.

    \appendix
    
    \section{Pseudo-code of the BLAS-implementation}
    Here, we describe the main integration routine presented in the main text. We do not include the Magnus expansion as this can be seen as adding effective control Hamiltonians and amplitudes to the problem. Those can be obtained according to eq.~\eqref{eq:magnus_coeffs}. The goal is to provide an easy way to port the integration parallelization to other platforms that offer Batched BLAS routines.
\begin{lstlisting}[caption=Implementation of the central integration routine in a C-style pseudo code]
// Working arrays X, D0, D1, length dim*dim*pts
// J[k] array with Bessel coefficients

// EXPANSION STEP
GEMM(OP_N, OP_N, 
dim*dim, pts, 1, 
1, &H0, dim*dim, 
c0, 1, 
0, &X, dim*dim);
GEMM(OP_N, OP_N, 
dim*dim, pts, amps, 
1, &H1, dim*dim, 
ck, pts, 
1, &X, dim*dim);

// PROPAGATE STEP
for (int k = MMAX; k >= 0; k--) {
    // D0 = D0 + 2 X @ D1 * dt
    GEMMStridedBatched(OP_N, OP_N,
    dim, dim, dim,
    dt,
    &X, dim, dim*dim,
    &D1, dim, dim*dim,
    -1,
    &D0, dim, dim*dim,
    pts);
    
    diagonal_add(D0, J[k], pts, dim);
    k--;
    if (k == 0) {
        ptr_accumulate = &handle->mtwo;
    }
    else {
        ptr_accumulate = &handle->mone;
    }
    GEMMStridedBatched(OP_N, OP_N,
    dim, dim, dim,
    dt,
    &X, dim, dim*dim,
    &D0, dim, dim*dim,
    ptr_accumulate,
    D1, dim, dim*dim,
    pts);
    diagonal_add(D0, J[k], pts, dim);
}
// D1 contains now the matrix exponentials

// REDUCTION STEP
complex *read = &D1;
complex *write = &D0;
complex *temp;

int remain_pts = pts;
int pad = 0;
while (remain_pts > 1) {
    pad = remain_pts %
    remain_pts = remain_pts/2;
    GEMMStridedBatched(OP_N, OP_N,
    dim, dim, dim,
    1, read, dim, dim*dim*2,
    read + dim*dim, dim, dim*dim*2,
    0, write, dim, dim*dim,
    remain_pts);
    if (pad > 0) {
        // One left over, need to copy to Array
        COPY(dim*dim,
        read + dim*dim*(remain_pts*2),
        1, write + dim*dim*(remain_pts), 1);
        remain_pts += 1;
    }
    temp = write;
    write = read;
    read = temp;
}
&D1 = read;
// D1 contains as first matrix the propagator

\end{lstlisting}

    \bibliographystyle{elsarticle-num}
    \bibliography{lit}

    \section*{Competing interests}
    The authors declare that they have no known competing financial interests or personal relationships that could have appeared to influence the work reported in this paper.
    
    \section*{Author contributions}
    K.H. initiated the project and implemented a prototype version of PARAMENT. K.H and P.W. implemented the production algorithm and cowrote the paper. All authors discussed the results.

    \section*{Acknowledgment}
    The authors thank Prof. Christian Degen for his support and proof-reading the manuscript. The authors acknowledge Prof. Oded Zilberberg and Timo Schmetzer for the fruitful discussions on the mathematical details of the numerics. 
    
    This work has been supported by Swiss National Science Foundation 
    (SNFS) Project Grant No. 200020 175600, the National Center of 
    Competence in Research in Quantum Science and Technology (NCCR QSIT), 
    and the Advancing Science and TEchnology thRough dIamond Quantum Sensing 
    (ASTERQIS) program, Grant No. 820394, of the European Commission.

\end{document}